\documentclass{article}
\usepackage{amsmath}
\usepackage{chicago}
\usepackage{graphicx}

\newtheorem{theorem}{Theorem}[section]

\makeatletter
\@addtoreset{equation}{section}
\makeatother

\newtheorem{example}{Example}[section]

\input{tcilatex}

\begin{document}

\title{On Bond Portfolio Management}
\author{Vladislav Kargin \thanks{%
Cornerstone Research, 599 Lexington Avenue, New York, NY 10022, USA;
slava@bu.edu} \thanks{%
The views expressed within this article are the views of the author and not
the views of Cornerstone Research.}}
\maketitle

\begin{abstract}
This paper describes a new method of bond portfolio optimization based on
stochastic string models of correlation structure in bond returns. The paper
shows how to approximate correlation function of bond returns, compute the
optimal portfolio allocation using Wiener-Hopf factorization, and check
whether a collection of bonds presents arbitrage opportunities.

\end{abstract}

\textbf{\newpage }

\section{Introduction}

A plot of monthly stock returns typically looks like the surface of the
stormy sea: Although stocks move up and down together with the market as the
sea level moves up and down reacting to the Moon's gravitation, the relative
movements of nearby points are independent. In contrast, returns on bonds
plotted against maturity look like the sea on a breezy but calm day -- the
nearby points are close to each other and move in accord like a wave. This
distinction is at the basis of the specificity of bond portfolio management:
Bonds are more structured but also more difficult to diversify. The present
paper contends that the best way to find the optimal bond portfolio is by
approximating the correlation structure in bond returns with a rational
function of the difference in maturities and reducing the optimization
problem to the inversion of an operator in a Hilbert space.

The importance of bond portfolio management is difficult to overestimate. In
2000 the values of government and corporate debt outstanding were 7.7 and 5
trillions respectively, compared with 17 trillions of corporate equity
outstanding. In addition, the relative importance of debt is rising: the
equity value in 2000 decreased by 13 percent from its 1999 level, while the
government and corporate debts were up 2 and 9 percent respectively.%
\footnote{%
Source: Securities Industry Association Factbook 2001, page 22.} Despite all
the effort to balance its budget, the federal government is still spending
more than it receives in revenues and therefore the likelihood that the
Treasury securities market will be shriveling is scant. So, given the
importance of debt markets how should investors optimize their bond
portfolios?

In the early 1980s \citeANP{heaney_cheng84} applied to bonds the techniques
for stock portfolio optimization. This approach, however, was soon abandoned
because bonds move together to a much greater extent than stocks and
modeling of these co-movements is harder.

Another idea appeared even earlier and got much greater application -- the
idea of immunization (\shortciteN{fisher_weil71}). The immunization
technique minimizes sensitivity of the portfolio with respect to small,
parallel shifts in all interest rates. So, this approach directly takes into
account the observation that interest rates are highly correlated. The
modern development of this idea uses stochastic programming (%
\shortciteN{dembo93}, \shortciteN{mulvey_zenios94}, \shortciteN{zenios95}, %
\shortciteN{zenios98}, \shortciteN{consigli_dempster98}, %
\shortciteN{zenios99}, \shortciteN{dupacova_bertocchi01}), in which the
investor formulates a set of scenarios for interest rate movements,
prescribes them probabilities, and minimizes a certain loss function -- for
example, loss that can occur with 5\% probability.

While practical, the immunization technique neglects certain aspects in the
dynamics of interest rates. The probabilities of scenarios are usually
extracted from a finite factor model, which does not capture all the
information about the statistical correlations available from the data. In
addition, the market models often used for generation of the scenario
probabilities are internally inconsistent since they require continual
re-calibration of parameters. What is needed is a better method for modeling
interest rate correlations.

The present paper explores a synthesis of the two methods of portfolio
optimization: The correlations are directly estimated using a plausible
assumption on their structure similar to the intuitive assumption of the
immunization technique. The estimated correlations are then used by the
standard portfolio optimization model.

This new method has its provenance in the random field models of bond
returns (\shortciteN{kennedy94}, \shortciteN{kennedy97}, %
\shortciteN{goldstein00}, and \shortciteN{santaclara_sornette01}), which
provide a more flexible framework than finite-factor models for estimation
correlations in bond returns. According to philosophy of these papers,
correlations in bond returns should be approximated as a smooth function of
the difference in maturities. Once the function is estimated, the problem of
optimization can be solved by inverting a special operator in a Hilbert
space spanned by bond returns, a task that was already extensively studied
in the communications engineering literature.

After the solution of the optimization problem is found, it can be
profitably applied to evaluate\ the opportunities that a given structure of
interest rates presents. A classical arbitrage opportunity arises when the
investor is willing to invest infinite amount in a security and his utility
from the investment is not bounded. The paper explains how to check for
existence of these opportunities. Since they rarely happen in the market,
the paper also explains how to check for the existence of near-arbitrage
opportunities, which are situations where the utility from the investment is
not infinite but abnormally large. It will be shown that there is a norm in
Hilbert space defined through parameters of the correlation structure such
that near-arbitrage opportunities arise only if the vector of expected
returns has large size in terms of this norm.

The rest of the paper is organized as follows. Section 2 explains
assumptions and notation. Section 3 is about Pad\'{e} approximations to the
correlation function. Section 4 solves the portfolio optimization in terms
of the Wiener-Hopf factorization. Section 5 is about arbitrage
opportunities. Section 6 computes the correlations and the optimal portfolio
for the Treasury interest rates in a certain period. And Section 7 concludes.

\section{Assumptions and Notation}

Recent research in the theory of interest rates treats the bond returns as a
random field that has a two-dimensional correlation structure. As usual,
bond returns follow a certain stochastic time process, but in addition the
contemporaneous returns of the bonds with different maturities are also
stochastically correlated. Most importantly, the contemporaneous
correlations are functions of maturities. I will assume that the
contemporaneous correlation of two bond returns is a smooth function of the
difference of their maturities. Let me introduce some notation to formulate
this assumption formally. Symbol $t$ will denote the time to maturity and
symbol $s$ the calendar time. Let $R(t,s)$ be a one-period return on the
bond with maturity $t$ at time $s:$%
\begin{equation}
R(t,s)=\log \frac{P(t-\Delta t,s+\Delta t)}{P(t,s)},
\end{equation}%
where $P(t,s)$ is the price of the bond with maturity $t$ at time $\ s$, and 
$\Delta t$ is the minimal possible difference between maturities. For
example, a month is a realistic minimal difference between maturities of
government bonds. Then the assumption claims that the correlation between
returns $R(t,s)$ and $R(t+\tau ,s+\varsigma )$ can be written as $C(\tau
,\varsigma ),$ where $\tau $ denotes difference between maturities and $%
\varsigma $ is difference between calendar times.

The assumption may be motivated by analogy with assumption of stationarity
in time series where autocovariance of the returns depends only on the
difference between times of these returns. On a deeper level, the market
percepts the bonds with close maturities as similar and the assumption says
that the relevant measure of similarity is the difference between
maturities. Another potentially useful measure of similarity for
non-government bonds would be the difference in credit ratings. Since the
question about appropriate measure of similarity is a question about market
perceptions, it needs further empirical investigation. A priori, the
assumption that is taken in this paper seems to be reasonable.

Let me in the following suppress the second argument in $R(t,s)$ if the
calendar time is unimportant, and in $C(\tau ,\varsigma )$ if the
correlation of contemporaneous returns is considered, $C(\tau )=:C(\tau ,0).$
Also, let $\Delta t=1.$ This notation is less cumbersome and no confusion
should arise.

It is useful to introduce artificial securities $S_{t}$ that have unit
variance and the following expected return 
\begin{equation}
E(t)=ER(t)/\sqrt{V(t)},
\end{equation}%
where $ER(t)$ is the expected return of bond with maturity $t$, and $V(t)$
is the variance of this return.\ Let also%
\begin{equation}
Y(t)=X(t)\sqrt{V(t)},  \label{Y_X}
\end{equation}%
where $X(t)$ is the holding in the bond with maturity $t.$ Then the bond
portfolio that holds $X(t)$ in the bond with maturity $t$ has the same
variance and expected return as the portfolio that holds $Y(t)$ in the
security $S_{t}.$ The usefulness of this reformulation is that securities $%
S_{t}$ have unit variance of return, and the optimization depends solely on
the correlation structure of the returns. Since it is the well-developed
correlation structure that makes the bond portfolio management different
from the stock portfolio management, the formulation it terms of normalized
holdings $Y(t)$ is useful as emphasizing the importance of the correlation
structure. When the optimal portfolio of securities $S_{t}$ is found, it is
easy to translate it back into the optimal portfolio of bonds by inverting
relationship (\ref{Y_X}).

The first step to the formulation of the investor optimization problem is to
introduce generating functions for correlations, holdings and expected
returns on securities $S_{t}:$%
\begin{eqnarray}
\widehat{C}(z) &=&:\sum_{\tau =0}^{\infty }C(\tau )z^{\tau }, \\
\widehat{Y}(z) &=&:\sum_{t=0}^{\infty }Y(t)z^{t},  \notag \\
\widehat{E}(z) &=&:\sum_{t=0}^{\infty }E(t)z^{t}.  \notag
\end{eqnarray}

The second step is to introduce some machinery of Hilbert spaces. The space
is needed because the number of bonds with different maturities is
potentially infinite. Let $\mathcal{H}$ be the linear space of the formal
series in variable $z$ with the bounded sum of squared coefficients:%
\begin{equation}
a(z)=\sum_{-\infty }^{\infty }a_{k}z^{k}\text{ such that }\sum_{-\infty
}^{\infty }|a_{k}|^{2}<\infty .
\end{equation}%
Scalar product $<a|b>=:\sum a_{i}\overline{b_{i}}$ turns $\mathcal{H}$ into
a Hilbert space. This scalar product can also be written in an integral form
often useful in computations:%
\begin{equation}
<a|b>=\frac{1}{2\pi i}\int_{|z|=1}a(z)\overline{b}(z^{-1})\frac{dz}{z}.
\end{equation}

Let $\mathcal{H}_{+}$ be a subspace of series with non-negative
coefficients, $\mathcal{H}_{-}$ its orthogonal complement, and $P_{+}$ the
orthogonal projector on $\mathcal{H}_{+}.$ Clearly, $\widehat{C}(z),\widehat{%
Y}(z),$ and $\widehat{E}(z)$ are elements of $\mathcal{H}_{+}.$ In addition,
the variance of a bond portfolio can be seen as a norm on $\mathcal{H}_{+}$.
To write down this norm explicitly and relate it to the function $\widehat{C}%
(z),$ we need to introduce multiplication operators: To any function 
\begin{equation*}
F(z)=\sum_{-\infty }^{\infty }f_{i}z^{i}\text{ such that }\sup
|f_{i}|<\infty ,
\end{equation*}%
corresponds an operator of multiplication by this function,%
\begin{equation*}
a(z)\rightarrow F(z)a(z).
\end{equation*}%
To make notation for multiplication operators distinct from notation for
corresponding functions, the operators will have a multiplication sign in
the subscript: $F_{\times }.$ So, function $F$ maps complex numbers to
complex numbers, and operator $F_{\times }$ maps the Hilbert space $\mathcal{%
H}$ to itself.

Using this notation it is easy to write the variance of the portfolio $Y(t):$%
\begin{equation}
\mathrm{Var}(Y)=<\widehat{Y}|P_{+}A_{\times }\widehat{Y}>,
\end{equation}%
where 
\begin{equation}
A(z)=\widehat{C}(z^{-1})+\widehat{C}(z)-1.
\end{equation}

Assume that the investor optimizes the short-run performance of his
portfolio. The myopia assumption allows to neglect correlation of bond
returns across time and reduce the problem to the static case. The dynamic
problem with non-myopic investor is considerably more difficult and deserves
further investigation.

Formally, the investor aims to maximize a certain linear combination of the
expected return and the variance in the return. We can write it as follows:%
\begin{equation}
U(Y)=<\widehat{E}|\widehat{Y}>-\gamma <\widehat{Y}|P_{+}A_{\times }\widehat{Y%
}>,
\end{equation}%
where $\gamma $ is a coefficient of risk aversion. This is the traditional
Markowitz portfolio optimization problem, written in the language of Hilbert
spaces. Note that it is precisely because of the assumption that the
correlations depend only on the difference in bond maturities that the
variance of the portfolio can be written as a scalar product of the
portfolio vector and its operator image.

The solution to this problem is:%
\begin{eqnarray}
\widehat{Y} &=&\frac{1}{2\gamma }\left[ P_{+}A_{\times }\right] ^{-1}%
\widehat{E},  \label{optimal_YU} \\
U &=&\frac{1}{4\gamma }<\widehat{E}|\left[ P_{+}A_{\times }\right] ^{-1}%
\widehat{E}>.  \notag
\end{eqnarray}%
It is important to note that similar solutions arise also for optimization
problems different from the problem considered here. For example, solutions
have the similar form for the problem of minimization of portfolio variance
with certain constraints on the portfolio composition. The problem in this
paper was chosen as the most familiar representative of this class of
problems but the method of solution is useful for all of them.

In all these problems the solution requires estimation of correlation
function $C(\tau )$ and inversion of the operator $P_{+}A_{\times }.$ We
will address the problem of inversion later and first explain how to
estimate the correlation function.

\section{Pad\'{e} Approximations}

Estimation of the correlation function is not a trivial task because the
data about returns of the bond with particular maturity are scarce.
Therefore, the data must be fitted by a smooth function. One way to approach
the fitting problem is to use Pad\'{e} approximations -- ratios of
polynomials constrained to have specific initial terms in their Taylor
expansions. Using the appropriate degrees of polynomials in numerator and
denominator, it is also possible to match any conjectured asymptotic
behavior of the function.

Formally, let $P_{M}(z)$ and $Q_{N}(z)$ be a couple of polynomials with the
ratio that have the Tailor expansion $f(z):$%
\begin{equation}
\frac{P_{M}(z)}{Q_{N}(z)}=\sum_{i=0}^{\infty }f(i)z^{i}.
\end{equation}%
This couple of polynomials is a Pad\'{e} approximation to the correlation
function $\widehat{C}(z)$ if the first $N+M+1$ coefficients of the Tailor
expansion coincide with the corresponding coefficients of $\widehat{C}(z):$%
\begin{equation}
f(i)=C(i)\text{ for }i=0,1,...,M+N.
\end{equation}%
Pad\'{e} approximations can be easily found by solving a system of linear
equations (see \shortciteN{baker_graves-morris96} for more information about
Pad\'{e} approximations).

A generalization of Pad\'{e} approximations may be helpful when the data is
noisy. The generalization requires that the coefficients of Taylor
expansions be only approximately equal:%
\begin{equation}
f(i)=\widehat{C}(i)+\varepsilon _{i}\text{ for }i=0,1,...,M+N+K.
\label{general_Pade}
\end{equation}%
By definition, the $[M,N,K]-$order generalized Pad\'{e} approximation
minimizes the sum of squared errors in (\ref{general_Pade}). The benefit of
this generalization is that it allows using correlations of bonds with
larger differences in maturities for more precise estimation of the
coefficients of the approximation.

\section{Wiener-Hopf factorization}

In the engineering literature the operator $P_{+}A_{\times }$ is called
rational filter, and its inversion is a well-known problem. It can be solve
by several efficient methods (\shortciteN{kailath_sayed_hassibi00}), from
which one of the most elegant is given by the Wiener-Hopf factorization. Let 
$\ln A(z)$ be decomposed as follows: 
\begin{equation}
\ln A(z)=A_{+}(z)+A_{-}(z),\text{ where }a_{+}(z)\in \mathcal{H}_{+}\text{
and }a_{-}(z)\in \mathcal{H}_{-.}
\end{equation}%
Then the Wiener-Hopf factorization theorem claims that 
\begin{equation}
\left[ P_{+}A_{\times }\right] ^{-1}=\left[ \exp (-A_{+})\right] _{\times
}P_{+}\left[ \exp (-A_{-})\right] _{\times }.
\end{equation}%
See \shortciteN{lax02} for the proof and Appendix \ref%
{appendix_factorization} for computational details. The benefit of this
theorem is that it allows explicitly inverting the infinite matrix
corresponding to operator $P_{+}A_{\times },$ which greatly reduces demand
for necessary computational resources.

The technique of the Wiener-Hopf factorization allows writing analytic
expressions for the optimal bond portfolio allocation and corresponding
utility:

\begin{theorem}
\label{optimal_portfolio}The optimal allocation is 
\begin{equation*}
\widehat{Y}=\frac{1}{2\gamma }\left[ \exp (-A_{+)}\right] _{\times }P_{+}%
\left[ \exp (-A_{-})\right] _{\times }\widehat{E}
\end{equation*}%
The corresponding utility function is 
\begin{equation*}
U=\frac{1}{4\gamma }<\widehat{E}|\left[ \exp (-A_{+})\right] _{\times }P_{+}%
\left[ \exp (-A_{-})\right] _{\times }\widehat{E}>.
\end{equation*}
\end{theorem}

\textbf{\noindent Proof:} This theorem is a direct consequence of the
Wiener-Hopf factorization theorem and expressions for optimal portfolio and
utility (\ref{optimal_YU}).

\begin{example}
AR(1) correlations and expectations
\end{example}

\noindent Let correlations\ between bond returns be as they are in AR(1)
time series model:%
\begin{equation}
\widehat{C}(z)=1+\sum_{i=1}^{\infty }\alpha ^{i}z^{i}=\frac{1}{1-\alpha z}.
\end{equation}%
Then 
\begin{eqnarray}
A(z) &=&\widehat{C}(z^{-1})+\widehat{C}(z)-1=\frac{1-\alpha ^{2}}{(1-\alpha
z)(1-\alpha z^{-1})}, \\
A_{+} &=&\ln \frac{1-\alpha ^{2}}{1-\alpha z},\text{ and }A_{-}=\ln \frac{1}{%
1-\alpha z^{-1}}.  \notag
\end{eqnarray}%
Therefore, 
\begin{equation}
\left[ P_{+}A_{\times }\right] ^{-1}=\frac{1}{1-\alpha ^{2}}(1-\alpha
z)_{\times }P_{+}(1-\alpha z^{-1})_{\times }.
\end{equation}%
Assume also for the purposes of this example that the normalized
expectations for bonds with longer maturities are smaller -- perhaps because
of large variance of the returns on longer maturity bonds. More precisely,
let the normalized expectations decline exponentially:%
\begin{equation}
\widehat{E}(z)=E_{0}(1+\sum_{i=1}^{\infty }\beta ^{i}z^{i})=\frac{E_{0}}{%
1-\beta z},
\end{equation}%
where $\beta <1$ is the rate of decline. Then, according to Theorem \ref%
{optimal_portfolio},

\begin{equation}
\widehat{Y}=\frac{E_{0}}{2\gamma }\frac{1-\alpha \beta }{1-\alpha ^{2}}\frac{%
1-\alpha z}{1-\beta z}\text{,}  \label{example_Y}
\end{equation}%
and%
\begin{eqnarray}
U &=&\frac{E_{0}^{2}}{4\gamma }\frac{1-\alpha \beta }{1-\alpha ^{2}}<\frac{1%
}{1-\beta z}|\frac{1-\alpha z}{1-\beta z}> \\
&=&\frac{E_{0}^{2}}{4\gamma }\frac{1-\alpha \beta }{1-\alpha ^{2}}\frac{1}{%
2\pi i}\int_{|z|=1}\frac{1}{1-\beta z^{-1}}\frac{1-\alpha z}{1-\beta z}\frac{%
dz}{z}  \notag \\
&=&\frac{E_{0}^{2}}{4\gamma }\frac{(1-\alpha \beta )^{2}}{(1-\alpha
^{2})(1-\beta ^{2})}.  \notag
\end{eqnarray}%
Note the wonderful symmetry of the expression relative to parameters that
govern expectations and correlations of bond returns. The symmetry
illustrates the idea that the investor will take into account both the
expectations and correlations of future returns.

Examination of the expression for the optimal portfolio allocation (\ref%
{example_Y}) reveals that the investor will sell short all bonds except the
bond with the shortest maturity provided that $\beta <\alpha .$ In this case
the bonds with larger maturities are valid only as hedging instruments. On
the contrary, he will buy the bonds with all the maturities if $\beta
>\alpha .$

\section{Arbitrage opportunities}

Classically, the absence of arbitrage opportunities means that there is no
investment with zero risk and positive return. In terms of the Hilbert space
formalism, zero-risk investments lie in the kernel of the operator $%
P_{+}A_{\times },$ and investments with zero expected return are orthogonal
to the vector of normalized expectations $\widehat{E}.$ Therefore, the
no-arbitrage condition on the structure of interest rates is 
\begin{equation}
\ker P_{+}A_{\times }\perp \widehat{E}.
\end{equation}

The classical arbitrage, however, can exist only if there is a perfect
correlation between a bond and a linear combination of other bonds, which is
unlikely to happen in reality. A formal tool to check for existence of
perfect correlations is given by a theorem of Hartman-Wintner-Widom (%
\shortciteN{hartman_wintner54}, and \shortciteN{widom60}), according to
which the operator $P_{+}A_{\times }$ is invertible if and only if $0\notin
\lbrack \mathrm{ess}\inf_{|z|=1}A(z),\mathrm{ess}\sup_{|z|=1}A(z)].$ In the
likely event that the operator $P_{+}A_{\times }$ is invertible, the
classical arbitrage opportunities are absent.

However, it seems natural to rule out also the near-arbitrage opportunities,
which are interest rates structure that permit to get utility greater than a
certain ``normality'' threshold. The near-arbitrage structures would allow
either getting a moderate excess return with very small risk, or a large
excess return with moderate risk. The formal criterion for deciding whether
an interest rate structure presents a near-arbitrage opportunity is
presented in the following Theorem:

\begin{theorem}
The market of bonds do not present near arbitrage opportunities if and only
if%
\begin{equation*}
<\widehat{E}|\left[ \exp (-A_{+})\right] _{\times }P_{+}\left[ \exp (-A_{-})%
\right] _{\times }\widehat{E}>\leq C,
\end{equation*}%
where C is a constant that depends on risk aversion of a typical investor.
\end{theorem}

\textbf{\noindent Proof: }This theorem directly follows from the definition
of near-arbitrage opportunity, and the expression for utility of optimal
portfolio.

\section{Application}

We use Treasury interest rates data by J. Huston McCulloch that represent
the 67 months from 8/1985 to 2/1991. These data give the zero-coupon yield
curve implicit in coupon bond prices. The yields have been defined for each
month using interpolation by cubic splines. From these data the returns on
holding a particular bond for one month have been computed.

The correlations have been estimated according to the formula%
\begin{equation}
C(\tau )=\frac{1}{N(s)N(t)}\sum_{s,t}(E_{s,t}-\overline{E_{s,t}}%
)(E_{s,t+\tau }-\overline{E_{s,t+\tau }}),
\end{equation}%
where $E_{s,t}$ are returns normalized by their standard deviation, and $%
N(s) $ and $N(t)$ are number of dates and maturities available for
estimation.

Figure 1 shows correlations predicted by a classical Pad\'{e} approximation
and the actual estimates of the correlations. Figure 2 shows correlations
from a generalized $\mathrm{Pad\acute{e}}$ approximation. From the
comparison of these figures, it is clear that the classical approximation is
good for small differences in maturities but severely underestimate the
correlation between bonds with larger difference in maturities. The
generalized $\mathrm{Pad\acute{e}}$ approximation is more balanced in the
sense that it approximates equally well the correlations for all differences
in maturities. On the other hand, the generalized approximation
underestimate the correlations between bonds with small difference in
maturities.

As usual in finance, evaluation of expected returns is more tricky than
estimation of covariances. In particular, it depends on what theory the
researcher holds about formation of interest rates. The results would be
different for segmentation, liquidity preference and expectation theories.
One possibility, assumed here for the purposes of illustration, is that
expectations of the future interest rate curve coincide with the current
interest rate curve. This allows to estimate expected return as follows:%
\begin{equation}
ER(t,s)=\log \frac{P(t-1,s)}{P(t,s)},
\end{equation}%
where $P(t,s)$ is the price of the bond with maturity $t$ at time$\ s$. It
should be emphasized that this is only a possible choice among many others.
It is appropriate for illustrative purposes because of its simplicity.

Figure 3 compares results of investment in optimal portfolio calculated
using a generalized Pade approximation with a benchmark. The benchmark is
results of the portfolio calculated under assumption that the bond returns
are uncorrelated. In the calculation of the portfolios, instead of
introducing a specific risk aversion parameter, the sum of investments is
constrained to be 1. Intuitively, this assumption means that even if the
investors buys more than his capital the excess is offset by a short sale of
a similar security. So, the return on such a portfolio is a return on the
unleveraged capital.

From Figure 3 it is clear that the optimal portfolio performs much better
than the benchmark portfolio. It has much lower variance and its monthly
returns are always positive in a striking difference with the returns of
benchmark portfolio. These results suggest that modelling the correlation
structure pays off.

\section{Conclusion}

An investor entering the business of bond portfolio management faces the
business that was regulated as early as in the time of Hammurabi when the
law required putting to death as a thief any man who received a deposit from
a minor or a slave without power of attorney, but that is still shaken by
demises of huge hedge funds, the business that attracts more bright
mathematicians and physicists than all mathematics and physics departments
in the country, that is operated by traders who play toy machine guns during
their lunch, and that demands deeper economic insight than stock portfolio
management will ever require, -- but when he faces this fascinating
business, the investor may perhaps be comforted by the thought that the
business is the most scientific and precise in the whole area of financial
speculation.

The present article is a contribution that shows how engineering techniques
can be applied to calculating optimal bond portfolios, and illustrates the
obtained formulas by developing a test for existence of arbitrage
opportunities. The method used in the article is a blend of the traditional
portfolio optimization with the techniques for estimation of correlations
between bonds of the similar maturities. While the present paper uses a
particular form of approximation, the method has more general applicability
and can be used with other estimation techniques.

\appendix

\section{Computation of Wiener-Hopf Factorization}

\label{appendix_factorization}Let$A(z)$ be a ratio of polynomials:%
\begin{equation*}
A(z)=a_{0}\frac{\prod_{i=1}^{M}(z-\theta _{i})}{\prod_{j=1}^{N}(z-\eta _{i})}%
.
\end{equation*}%
Let $\theta _{i}^{+}$ and $\eta _{i}^{+}$ be zeros and poles outside the
unit circle, and $\theta _{i}^{-}$ and $\eta _{i}^{-}$ be zeros and poles
inside the unit circle. Then 
\begin{eqnarray*}
A_{+}(z) &=&\ln \left[ a_{0}\frac{\prod (z-\theta _{i}^{+})}{\prod (z-\eta
_{i}^{+})}\right] , \\
A_{-}(z) &=&\ln \left[ \frac{\prod (z-\theta _{i}^{-})}{\prod (z-\eta
_{i}^{-})}\right] .
\end{eqnarray*}

\bibliographystyle{CHICAGO}
\bibliography{comtest}

\newpage 

\FRAME{ftbpFUO}{4.8179in}{3.8043in}{0pt}{\Qct{}\Qcb{Actual and Fitted
Correlation Functions for Classical Pad\'{e} Approximation}}{\Qlb{figure1}}{%
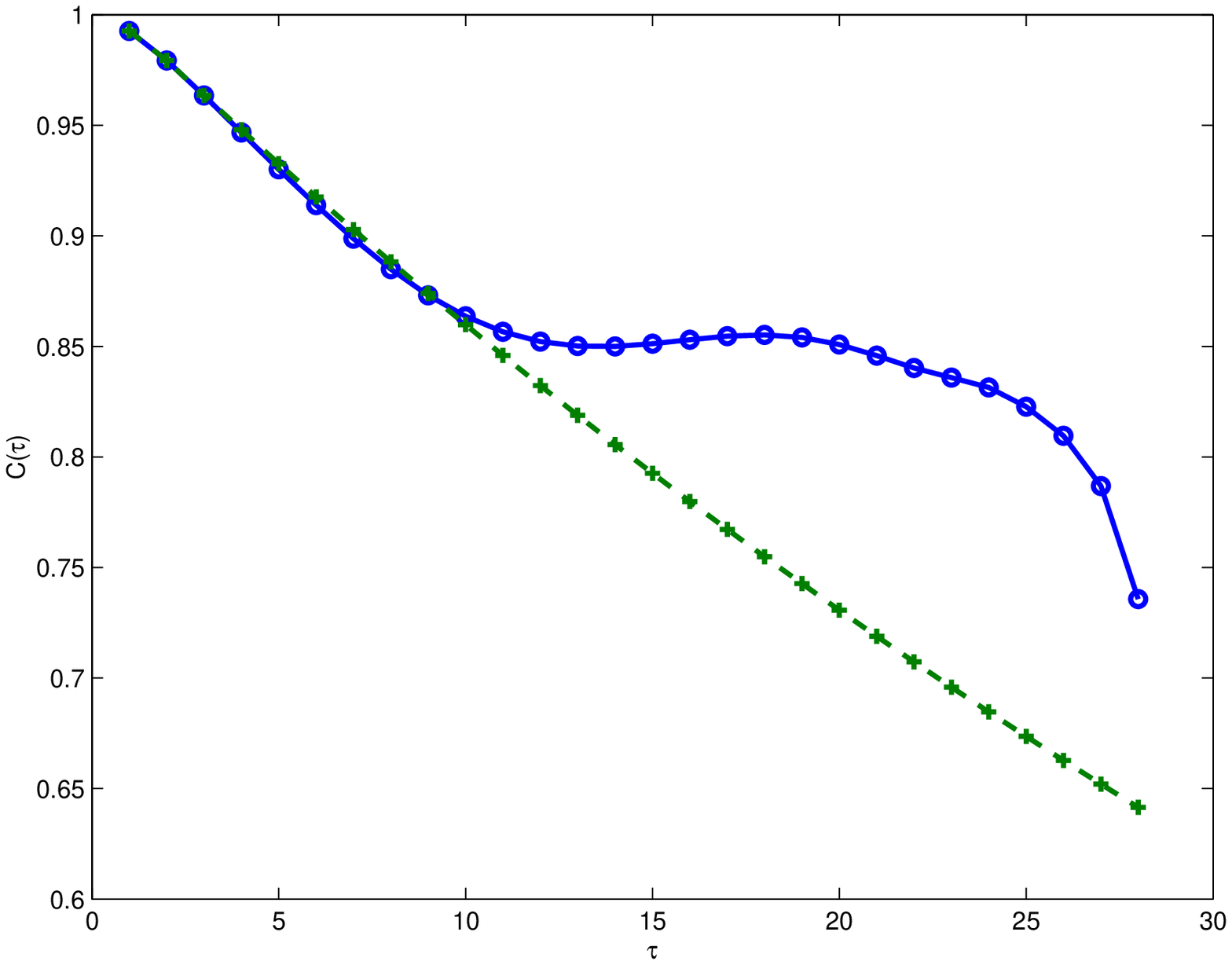}{\special{language "Scientific Word";type
"GRAPHIC";maintain-aspect-ratio TRUE;display "USEDEF";valid_file "F";width
4.8179in;height 3.8043in;depth 0pt;original-width 6.845in;original-height
5.3956in;cropleft "0";croptop "1";cropright "1";cropbottom "0";filename
'../Pictures for Papers/figure1_bonds_paper.eps';file-properties "XNPEU";}}

The solid line marked by circles is the estimate of the actual correlation
function for bonds of different maturities. The dashed line marked by pluses
is the fitted correlations from the classical Pad\'{e} [1/2] approximation.
The vertical axis is correlations. The horizontal axis is the differences
between maturity times measured in years. \newpage 

\FRAME{ftbpFUO}{4.8196in}{3.8251in}{0pt}{\Qct{}\Qcb{Actual and Fitted
Correlation Functions for Generalized Pad\'{e} Approximation}}{\Qlb{figure2}%
}{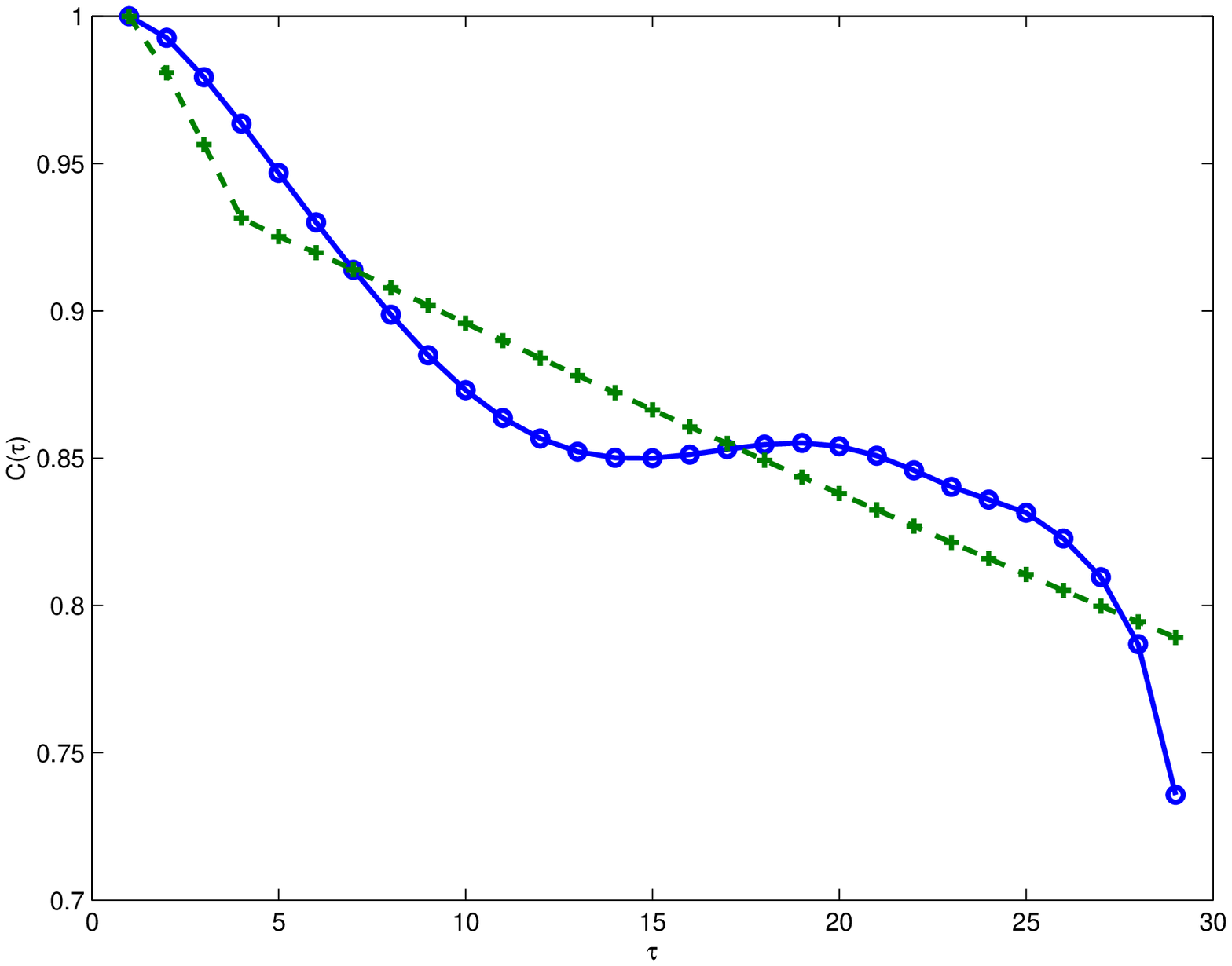}{\special{language "Scientific Word";type
"GRAPHIC";maintain-aspect-ratio TRUE;display "USEDEF";valid_file "F";width
4.8196in;height 3.8251in;depth 0pt;original-width 6.845in;original-height
5.4232in;cropleft "0";croptop "1";cropright "1";cropbottom "0";filename
'../Pictures for Papers/figure2_bonds_paper.eps';file-properties "XNPEU";}}

The solid line marked by circles is the estimate of the actual correlation
function for bonds of different maturities. The dashed line marked by pluses
is the fitted correlations from the generalized Pad\'{e} [0/5/28]
approximation. The vertical axis is correlations. The horizontal axis is the
differences between maturity times measured in years.

\newpage 

\FRAME{ftbpFUO}{4.868in}{3.7957in}{0pt}{\Qct{}\Qcb{Returns for Optimal and
Benchmark Portfolios}}{\Qlb{figure3}}{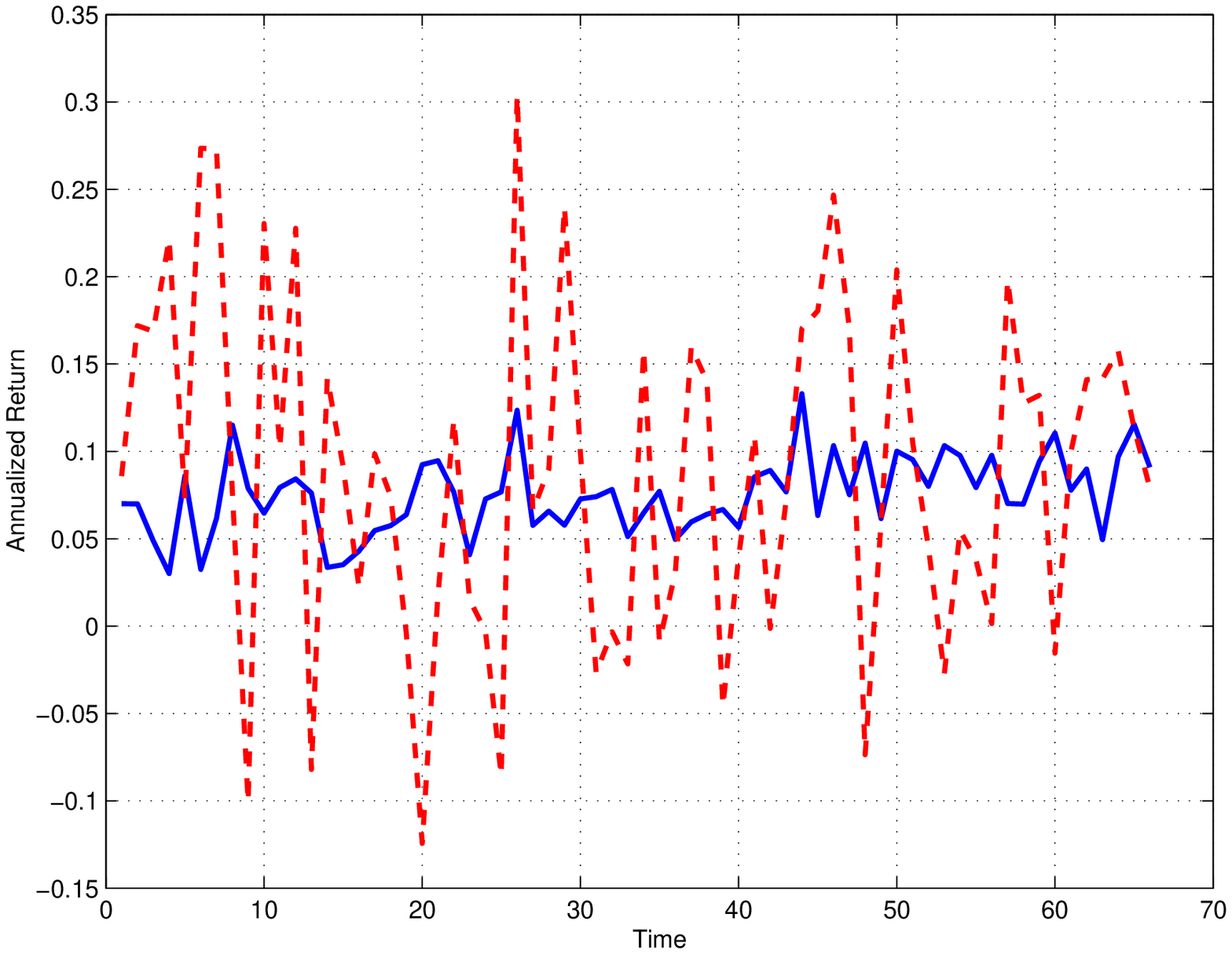}{\special%
{language "Scientific Word";type "GRAPHIC";maintain-aspect-ratio
TRUE;display "USEDEF";valid_file "F";width 4.868in;height 3.7957in;depth
0pt;original-width 6.915in;original-height 5.3835in;cropleft "0";croptop
"1";cropright "1";cropbottom "0";filename '../Pictures for
Papers/figure3_bonds_paper.eps';file-properties "XNPEU";}}

The solid line shows the returns of the optimal portfolio computed using
[0/5/28] generalized Pad\'{e} approximation. The dashed line is the returns
of the benchmark portfolio computed under assumption of zero correlations
between bond returns. The horizontal axes shows time in months; the vertical
axes shows annualized monthly returns.

\end{document}